\documentclass[12pt]{article}

\usepackage{amsmath,amssymb,amsthm}
\usepackage{epsfig}
\usepackage{arydshln}

\newcommand{\bfm}[1]{\mbox{\boldmath ${#1}$}}
\newcommand{\nonum}{\nonumber \\}
\newcommand\eq[1] {(\ref{#1})} 
\newcommand{\beqa}{\begin{eqnarray}}
\newcommand{\eeqa}[1]{\label{#1}\end{eqnarray}}
\newcommand{\beq}{\begin{equation}}
\newcommand{\eeq}[1]{\label{#1}\end{equation}}
\newcommand{\Grad}{\nabla}
\newcommand{\Div}{\nabla \cdot}
\newcommand{\Curl}{\nabla \times}

    \newcommand{\Imag}{\mathop{\rm Im}\nolimits}
\newcommand{\Tr}{\mathop{\rm Tr}\nolimits}

\newcommand{\lang}{\langle} 
\newcommand{\rang}{\rangle}
\newcommand{\Md}{\partial}

\newcommand{\Gd}{\delta}
\newcommand{\Ge}{\epsilon}
\newcommand{\Gve}{\varepsilon}

\newcommand{\Gc}{\chi}

\newcommand{\Gk}{\kappa}

\newcommand{\Gm}{\mu}

\newcommand{\Gt}{\theta}

\newcommand{\Gr}{\rho}
 
\newcommand{\Gs}{\sigma}

\newcommand{\BGe}{\bfm\epsilon}

\newcommand{\BGr}{\bfm\rho}

\newcommand{\BGs}{\bfm\sigma}

\newcommand{\CE}{{\cal E}}

\newcommand{\CJ}{{\cal J}}

\newcommand{\CM}{{\cal M}}

\newcommand{\CS}{{\cal S}}

\newcommand{\bpm}{\begin{pmatrix}}
\newcommand{\epm}{\end{pmatrix}}

\usepackage{url}
\usepackage{hyperref}

\def\Ba{{\bf a}}
\def\Bb{{\bf b}}

\def\Be{{\bf e}}

\def\Bh{{\bf h}}

\def\Bj{{\bf j}}

\def\Bs{{\bf s}}

\def\Bu{{\bf u}}
\def\Bv{{\bf v}}

\def\Bx{{\bf x}}

\def\BA{{\bf A}}
\def\BB{{\bf B}}
\def\BC{{\bf C}}

\def\BE{{\bf E}}
\def\BF{{\bf F}}

\def\BI{{\bf I}}
\def\BJ{{\bf J}}
\def\BK{{\bf K}}
\def\BL{{\bf L}}

\def\BN{{\bf N}} 

\def\BP{{\bf P}}
\def\BQ{{\bf Q}}
\def\BR{{\bf R}}

\def\BV{{\bf V}}
\def\BW{{\bf W}}

\usepackage{graphics}
\usepackage{amssymb}

\title{Some open problems in the theory of composites}
\author{Graeme W. Milton}
\date{\small{Department of Mathematics, University of Utah, Salt Lake City, UT 84112, USA
\\Email: milton@math.utah.edu}}
\begin{document}
\maketitle
\vspace{2ex}
\begin{abstract}
  A selection of open problems in the theory of composites is presented. Particular attention is drawn to
  the question of whether two-dimensional, two-phase, composites with general geometries have the same set of possible effective
  tensors as those of hierarchical laminates. Other questions involve the conductivity and elasticity of composites.
  Finally some future directions for wave and other equations are mentioned. 
\end{abstract}
\vspace{3ex}
\section{Introduction}
\setcounter{equation}{0}

The theory of composite materials has seen a resurgence of interest thanks to the discovery of novel properties
and a dramatic rise in our ability to manufacture desired microgeometries: see for instance the review \cite{Kadic:2019:M}
and references therein. Back in the 1980's and 1990's there was also a rapid increase in interest, partly due to the
recognition that the solution of optimal design problems often require composite microstructures in the
design. This gave rise to the area of topology optimization which has had enormous impact, moving into
the mainstream of engineering design: see, for example, the book \cite{Bendsoe:2004:RPM}. From a mathematics perspective there were accompanying
rapid developments: in our understanding of homogenization, which underlies the use of effective
moduli to describe macroscopic responses; in bounds on effective moduli, coupled with
the identification of microstructures that attain them; in the theory governing microgeometry independent
exact relations satisfied by effective moduli; and in the discovery of composites with unexpected
properties, as surveyed in the books
\cite{Bensoussan:1978:AAP, Zhikov:1994:HDO, Cherkaev:2000:VMS, Torquato:2001:RHM, Allaire:2000:SOH, Milton:2002:TOC, Tartar:2009:GTH, Milton:2016:ETC,
  Grabovsky:2016:CMM}.

Given the recent interest it is perhaps appropriate to draw attention to some of the open problems
generated in the mathematical research that is now mostly over three decades old, as well as questions
generated by more recent investigations. The problems here are by no means exhaustive. Rather they are ones I have encountered in my research work and
found quite difficult, usually because I have no idea how to solve them. Some are just of theoretical
interest, while others should be of interest to both experimentalists and theorists alike. 
The problems reflect my own research interests, both past and present, and other experts in the field would undoubtedly choose
a different set. Many are old outstanding problems, where it is difficult to dig in the hard soil, but some address new topics
where the soil is more fertile and it is easier to break ground.

\section{Open problems involving quasiconvexification}
\setcounter{equation}{0}

Here we present a selection of open problems that are related to quasiconvexification. For a recent survey
of selected results pertaining to quasiconvexity, and the closely related topic of weak lower semicontinuity, see \cite{Dacorogna:2007:DMC, Benesova:2017:WLS} and references therein.
The focus is largely on two-phase composites, and the corresponding two-well quasiconvexification problems, since these are perhaps of greatest interest
in the field of composites (though some effects, such as getting negative or unbounded thermal expansion coefficients
from materials having only positive thermal expansion coefficients,
require at least three phases \cite{Lakes:1996:CSS,Sigmund:1997:DME}). In this age of 3d-printing it is now relatively easy to manufacture tailored microstructures of one phase
plus void that can then be infilled
to obtain a two-phase material. One is interested in the range the effective tensors can have as the microgeometry varies over all
configurations. This range is known as the $G$-closure and provides limits for what one can expect to achieve when one tries to optimize the local
response using relatively simple practical microstructures obtained, for example, by topology optimization. The question we explore is whether it suffices to consider only
hierarchical laminate geometries rather all conceivable microstructures. Hierarchical laminate geometries have the advantage that it is relatively easy to calculate their
effective properties (see, for example, \cite{Tartar:1985:EFC, Francfort:1986:HOB}, Chapter 9 in \cite{Milton:2002:TOC}, and references therein).
\newline
\medskip

\noindent
We start with:
\newline
\medskip
\noindent
{\it Problem 1: Is rank convexity equal to quasiconvexity for the two well problem in two spatial dimensions?}
\newline
\medskip

\noindent
Given two self-adjoint positive definite mappings $\BL_1$ and $\BL_2$ on the space $\CS_{m}$ of real $2\times m$ matrices, equipped with the standard inner product
\beq \BA_1\cdot\BA_2=\Tr(\BA_1\BA_2 ^T), \eeq{0.1}
where $\Tr$ denotes the trace, and $\BA_1,\BA_2\in\CS_{m}$, and given $\BF_1,\BF_2\in\CS_{m}$, and two reals $c_1$ and $c_2$, consider
the two well ``energy'',
\beq W(\BF)=\min\{W_1(\BF),W_2(\BF)\},\quad \BF\in\CS_{m}, \eeq{1.1}
where the $W_j(\BF)$, $j=1,2$, are the quadratic wells
\beqa
W_j(\BF)& = & (\BF-\BF_j)\cdot\BL_j(\BF-\BF_j)+k_j \nonum
& = &\BF\cdot\BL_j\BF+2\BV_j\cdot\BF+c_j,\quad\BV_j=-\BL_j\BF_j,\quad c_j=k_j+\BF_j\cdot\BL_j\BF_j.
\nonum &~&
\eeqa{1.1a}
The quasiconvexification of $W(\BF)$ is given by
\beq QW(\BF)=\inf_\Bu\lang W(\BF+\Grad\Bu)\rang,
\eeq{1.2}
where the infimum is over all $m$-component periodic potentials $\Bu(\Bx)$ and the average $\lang\cdot\rang$ is over the unit cell of periodicity.
(We adopt the convention that the elements of $\Grad\Bu$ are $\{\Grad\Bu\}_{ij}=\Md u_j/\Md x_i$.)
An energy $W_0(\BF)$ is said to be rank-one convex if
\beq W_0(\Ba\otimes\Bb)\geq pW_0(\Ba\otimes\Bb)+(1-p)W_0(\Ba\otimes\Bb), \eeq{1.3}
for all real $p\in[0,1]$, all real $2$-component vectors $\Ba$, and all real $m$-component vectors $\Bb$. The rank-one convexification of $W(\BF)$, denoted $RW(\BF)$,
is the highest rank-one convex energy that lies equal or below $W(\BF)$ for all $\BF$. So the question is whether $QW(\BF)=RW(\BF)$
for all choices of $m, \BK_1,\BK_2,\BF_1,\BF_2,c_1,c_2$? We will see that this can be reduced to the problem with $\BF_1=\BF_2=0$. 
Clearly the problem does not change if we add the same constant to $c_1$ and $c_2$. So without loss of generality we can assume that $c_1$ and $c_2$ are sufficiently
large so that
\beq \BK_1=\bpm \BL_1 & \BV_1 \\  \BV_1^T & c_1 \epm>0,\quad \BK_2=\bpm \BL_2 & \BV_2 \\  \BV_2^T & c_2 \epm>0.
\eeq{1.4}
In terms of these we have
\beq W_j(\BF)=\bpm \BF \\ 1\epm\cdot\BK_j\bpm \BF \\ 1\epm, \eeq{1.5}
in which the inner product is the obvious generalization of \eq{0.1}.

In the field of composites problem 1 is equivalent to the following question:
\newline
\medskip

\noindent
{\it Problem 2: For two-phase composites in two spatial dimensions, such that phase 1 occupies a volume fraction $f$,
  is the $G_f$-closure equal to its lamination closure when the fields on the right of the constitutive law
  have $n$ components, each being the sum of a real $2$ component vector and the gradient of a scalar periodic potential, while the fields on the left of the constitutive law also
  have $n$ components, each having zero divergence, in which $n$ is an arbitrary positive integer?}
\newline
\medskip

\noindent
The constitutive law takes the form
\beq  \underbrace{\bpm \Bj^{(1)}(\Bx) \\ \Bj^{(2)}(\Bx) \\ \vdots \\ \Bj^{(k)}(\Bx)\epm}_{\BJ(\Bx)}=\BL(\Bx)
\underbrace{\bpm \Be^{(1)}(\Bx) \\ \Be^{(2)}(\Bx) \\ \vdots \\ \Be^{(k)}(\Bx) \epm}_{\BE(\Bx)},
\eeq{1.10}
where the $\Bj^{(i)}(\Bx)$, $\Be^{(j)}(\Bx)$, $\BL(\Bx)$ all have the same periodicity and satisfy
\beq \Div\Bj^{(i)}=0,\quad \Be^{(j)}=\Be^{(j)}_0+\Grad V_j,\quad
\BL(\Bx)=\Gc(\Bx)\BL_1+[1-\Gc(\Bx)]\BL_2, \eeq{1.11}
in which the $\Be^{(j)}_0$ are constant vectors, the $V_j(\Bx)$ are periodic potentials, $\Gc(\Bx)$ is the indicator function
\beqa \Gc(\Bx)& = & 1\quad\text{in phase 1}, \nonum
              & = & 0\quad\text{in phase 2},
\eeqa{1.12}
satisfying $\lang\Gc\rang=f$, in which the angular brackets $\lang \,\,\rang$ denote a volume average over the unit cell of periodicity, 
and $\BL_1$ and $\BL_2$ are self-adjoint positive definite mappings on $\CS_{n}$.
Thus $\BL_1$ and $\BL_2$ take the block matrix form
\beq \BL_j=\bpm \BGs^{(11)}_j & \BGs^{(12)}_j & \ldots & \BGs^{(1n)}_j \\
\BGs^{(21)}_j & \BGs^{(22)}_j & \ldots & \BGs^{(2n)}_j & \\
\vdots & \vdots & \ddots & \vdots & \\
\BGs^{(n1)}_j & \BGs^{(n2)}_j & \ldots & \BGs^{(nn)}_j 
\epm,\quad j=1,2,
\eeq{1.13}
in which each $\BGs^{(k\ell)}_j$ is a $2\times 2$ matrix, with $\BGs^{(k\ell)}_j=[\BGs^{(\ell k)}_j]^T$. The linear relation
\beq \lang\BJ\rang=\BL_*\lang\BE\rang \eeq{1.14}
determines the effective tensor $\BL_*$. The $G_f$-closure, $G_f(\BL_1,\BL_2)$, is the closure of the set of values $\BL_*$ takes
as $\Gc(\Bx)$ ranges over all possible indicator functions satisfying $\lang\Gc\rang=f$. In other words the microstructure
varies over all possible configurations in which phase 1 occupies a volume fraction $f$. The lamination closure, $G^L_f(\BL_1,\BL_2)$
is the closure of the set of values $\BL_*$ takes as $\Gc(\Bx)$ ranges over the indicator functions of multiple-rank laminate materials
satisfying $\lang\Gc\rang=f$. Multiple-rank laminate materials are hierarchical materials, obtained by an iterative process of lamination in different directions
on larger and larger length scales, ideally with an infinite ratio between the length scales at each stage of construction. A rank $1$ laminate
is just a simple laminate of the phases, which can be regarded as rank $0$ laminates. A rank $m$ laminate is obtained by layering
together a rank $m-1$ laminate with a laminate of rank $m-1$ or less. 
\newline
\medskip

\noindent
{\bf Remark 2.1}
\newline
\medskip

\noindent
The equivalence of $G_f(\BL_1,\BL_2)$ and  $G^L_f(\BL_1,\BL_2)$ in the case $n=1$ has been established by Nesi \cite{Nesi:1993:UQF} and
Grabovsky \cite{Grabovsky:1993:GCT,Grabovsky:1996:BEM}, subject to certain
assumptions about $\BL_1=\BGs^{(11)}_1$ and $\BL_2=\BGs^{(11)}_2$. (The $n=1$ case where $\BL_1$ and $\BL_2$ do not commute, and $\BL_1-\BL_2$ is neither positive
nor negative semidefinite, is unresolved to my knowledge). They built on earlier work of Lurie and Cherkaev \cite{Lurie:1982:AEC} and Murat and Tartar \cite{Murat:1985:CVH}
 who treated, using a variational approach known as the translation method, or method of compensated compactness, the
case where $\BGs^{(11)}_1$ and $\BGs^{(11)}_2$ are both proportional to the identity matrix, corresponding to isotropic materials.
For $n=2$ it is an open question as to whether they are equivalent.
In planar elasticity with two, possibly anisotropic, phases with fixed orientations, which is a subcase of the $n=2$ case,
existing evidence points to them being equivalent. In three-dimensional elasticity one needs microstructures, such as pentamode materials
\cite{Milton:1995:WET}, that are stiff with respect one loading, yet compliant with respect to all other loadings (which span a five-dimensional space),
and it is by no means clear that their behavior can be mimicked by hierarchical laminate structures. 
\newline
\medskip

\noindent
{\bf Remark 2.2}
\newline
\medskip

\noindent
In two spatial dimensions Grabovsky \cite{Grabovsky:2017:MIF} has an example of a manifold $\CM$ of tensors $\BL$ that is stable under lamination but not
under homogenization. This suggests that by picking anisotropic $\BL_1,\BL_2\in\CM$ one might find a $\Gc(\Bx)$ such that $\BL_*$ is not in $\CM$, thus establishing
that $G(\BL_1,\BL_2)$ and $G^L(\BL_1,\BL_2)$ differ. However, the analysis showing that $\CM$ is stable under lamination \cite{Grabovsky:1998:EREa} extends directly
to all two-phase composite geometries as can be seen from \cite{Grabovsky:2000:ERE} once takes the ``reference tensor'' $\BL_0$ equal to $\BL_2$.
We conclude that $\BL_*\in\CM$. The same analysis applies to
any manifold $\CM$ stable under lamination in any spatial dimension: if  $\BL_1,\BL_2\in\CM$ then also $\BL_*\in\CM$, for any indicator function $\Gc(\Bx)$, not just
those corresponding to laminate geometries. 
\newline
\medskip

\noindent
{\bf Remark 2.3}
\newline
\medskip

\noindent
If indeed $G(\BL_1,\BL_2)$ and $G^L(\BL_1,\BL_2)$ differ for some $n$ and some $\BL_1\geq 0$ and $\BL_2\geq 0$, the next questions become: can one identify
the minimum value $n_0$ of $n$ for which they differ for some $\BL_1$ and $\BL_2$, and given $n\geq n_0$ can one identify the set of pairs $(\BL_1,\BL_2)$ 
for which they differ, or for which $G_f(\BL_1,\BL_2)$ and $G^L_f(\BL_1,\BL_2)$ differ for fixed $f$? More generally, if one has a composite with $k$ phases,
what is the smallest value of $n$ for which $G(\BL_1,\BL_2,\ldots,\BL_k)$ and $G^L(\BL_1,\BL_2,\ldots,\BL_k)$ differ or for which $G(\BK_1,\BK_2,\ldots,\BK_k)$
and $G^L(\BK_1,\BK_2,\ldots,\BK_k)$ differ? A variant of an example of {\v{S}}ver{\'a}k \cite{Sverak:1992:ROC} shows that $G(\BK_1,\BK_2,\ldots,\BK_7)$
and $G^L(\BK_1,\BK_2,\ldots,\BK_7)$ differ when $n=3$ (see section 31.9 of \cite{Milton:2002:TOC}).
\newline
\medskip

\noindent
{\bf Remark 2.4}
\newline
\medskip

\noindent
In three spatial dimensions it seems quite likely that there are two phase geometries such that $G(\BL_1,\BL_2)$ and $G^L(\BL_1,\BL_2)$ differ.
 To obtain a candidate example, one considers the conductivity
equations in the presence of a small magnetic field $\Bh=(h_1,h_2,h_3)$. In a two-phase medium where phase 1 is isotropic while phase 2 is void, these take the form:
\beq \Bj(\Bx)=\Gc(\Bx)\BGr^{-1}\Be(\Bx),\quad\Div\Bj=0,\quad \Be=\Be_0+\Grad V,\quad\BGr=\Gr\BI+R^H\bpm 0 & -h_3 & h_2\\ h_3 & 0 & -h_1 \\ -h_2 & h_1 & 0 \epm, \eeq{1.15a}
where $R^H$ is the Hall coefficient of phase 1, and $\BGr$ is its resistivity tensor.
Assuming that the microstructure is isotropic or has cubic symmetry, the effective resistivity tensor $\BGr_*=\BGs_*^{-1}$ (if it exists) to first order in $\Bh$ takes the form
\beq \BGr_*=\Gr_*\BI+R^H_*\bpm 0 & -h_3 & h_2\\ h_3 & 0 & -h_1 \\ -h_2 & h_1 & 0 \epm. \eeq{1.15b}
Numerical results \cite{Kadic:2015:HES, Kern:2018:THE} and corresponding physical experiments \cite{Kern:2016:EES} show that in certain microstructures of interlinked tori, arranged to
have cubic symmetry,  $R^H_*$ and $R^H$ can have opposite signs. While it was commonly believed that the sign of Hall coefficient corresponds to the sign of the charge carrier,
these composites provide a counterexample
as they show the macroscopic  Hall coefficient can be opposite in sign to the Hall coefficients of the constituent materials, assuming their
Hall coefficients are zero or share a common sign. The argument that the Hall coefficient corresponds to the sign of the charge carrier assumes
that the electrons, or holes, travel in straight lines, which of course is not the case in these composite materials.
The microstructures were motivated
by a three-phase example \cite{Briane:2009:HTD} having cubic symmetry, where it was rigorously shown that the Hall coefficients $R^H_1, R^H_2$ and $R^H_3$ of all three isotropic phases can be non-negative, while at the same time
$R^H_*$ is negative. One can explain this \cite{Briane:2009:HTD,Kern:2018:THE} in terms of the ``matrix valued'' electric field $\BE(\Bx)$ whose three column vectors $\Be_1(\Bx), \Be_2(\Bx)$, and $\Be_3(\Bx)$ each
solve the conductivity equations, with zero magnetic field (i.e. the same $\Gc(\Bx)$ and $\BGr=\Gr\BI$).
Assuming $\lang\BE\rang=\BI$, a perturbation argument \cite{Bergman:1983:SDL, Briane:2009:HTD} shows
that the sign change of the Hall coefficient is related to the fact that the trace of the cofactor matrix of $\BE(\Bx)$ changes sign, at least in certain regions in
the unit cell of periodicity. On the other hand, in any multiple rank laminates (with $\lang\BE\rang=\BI$) Briane and Nesi show that the determinant of $\BE(\Bx)$ remains positive \cite{Briane:2004:WKL}, whereas it does take negative values in certain regions in the interlinked tori geometries \cite{Briane:2004:CSC}.
While they show that the trace of the cofactor matrix of $\BE(\Bx)$ can
change sign in three phase multiple rank laminates, it is an open question as to whether it can change sign in two phase multiple rank laminates. If it cannot, then the path is clear
to establishing that there are three-dimensional two phase geometries such that $G(\BL_1,\BL_2)$ and $G^L(\BL_1,\BL_2)$ differ.
We add that while in \eq{1.15a} the conductivity tensor $\BGs(\Bx)=\Gc(\Bx)\BGr^{-1}$ is not symmetric, one can perturb the problem slightly so that phase 2 is slightly
conducting, and then, using ideas of Cherkaev and Gibiansky \cite{Cherkaev:1994:VPC},
make a transformation to an equivalent problem where the tensor entering the constitutive law is real, symmetric, and positive definite
(see \cite{Milton:1990:CSP} and section 12.11 of \cite{Milton:2002:TOC}).
Also one can introduce a periodic vector potential $\Bv$ for $\Bj-\lang\Bj\rang$ in \eq{1.15a} so that $\Bj-\lang\Bj\rang$ is expressed in terms of the antisymmetric part of 
$\Grad \Bv$ using the completely
antisymmetric Levi-Civita tensor giving $\Bj-\lang\Bj\rang=\Curl\Bv$, while on the other hand the Levi-Civita tensor applied to $\Grad V$ gives an antisymmetric field
that has zero divergence. Then the equations can be manipulated into the same form as \eq{1.10}-\eq{1.13}  with
real $\BGs^{(k\ell)}_j=\BGs^{(\ell k)}_j$.
\newline
\medskip

\noindent
{\bf Equivalence between problems 1 and 2}
\newline
\medskip

\noindent
The connection between problems 1 and 2 is implicit in existing results. To see this,
we first consider a problem associated with, and in fact equivalent to, problem 2. This is to characterize the $G$-closure associated with the equations
\beq \bpm \BJ(\Bx) \\ \Bs(\Bx) \epm= \BK(\Bx)\bpm \BE(\Bx) \\ \Gt \epm,\quad \BK(\Bx)=[\Gc(\Bx)\BK_1+(1-\Gc(\Bx))\BK_2],
\eeq{1.15}
in which $\BJ(\Bx)$ and $\BE(\Bx)$ satisfy the same constraints as in problem 2, $\BK_1$ and $\BK_2$ are positive definite and given by \eq{1.4}, the indicator function
$\Gc(\Bx)$ is again given by \eq{1.12}, but not subject to the constraint that $\lang\Gc\rang=f$,
$\Gt$ is a constant scalar, and $\Bs(\Bx)$ is an arbitrary scalar valued function having the same periodicity as $\Gc(\Bx)$. The effective
tensor $\BK_*$ is defined by the linear relation
\beq  \bpm \lang\BJ\rang \\ \lang\Bs\rang \epm=\BK_*\bpm \lang\BE\rang \\ \Gt \epm.
\eeq{1.16}
The $G$-closure, $G(\BK_1,\BK_2)$ is the closure of the set of values $\BK_*$ takes as $\Gc(\Bx)$ ranges over all possible indicator functions, with no
constraint on the volume fraction.

Now when $\Gt=0$ \eq{1.15} when solved for $\BJ(\Bx)$ is exactly the same as \eq{1.10}. This implies that $\BK_*$ takes the form
\beq \BK_*=\bpm \BL_* & \BV_* \\  \BV_*^T & c_* \epm, \eeq{1.17}
where $\BL_*$ is the exactly the same effective tensor associated with problem 2, defined by \eq{1.14}. Furthermore if we assume
that $\BL_1-\BL_2$ is non-singular (by, if necessary, perturbing the problem) then we can find constant fields $\BJ(\Bx)=\BJ_0$
and $\BE(\Bx)=\BE_0$ that solve \eq{1.15}, and thus obtain formulas for $\BV_*$ and $c_*$. 
This is a standard technique in the theory of composites (see, for example, Chapter 5 and
in particular Section 5.4 in \cite{Milton:2002:TOC} and references therein). Specifically, \eq{1.15} and \eq{1.16} imply
\beqa \BJ_0& = & \BL_1\BE_0+\BV_1\Gt=\BL_2\BE_0+\BV_2\Gt=\BL_*\BE_0+\BV_*\Gt, \nonum
\lang\Bs\rang& = & [f\BV_1+(1-f)\BV_2]^T\BE_0+[fc_1+(1-f)c_2]\Gt=\BV_*\BE_0+c_*\Gt,\nonum &~&
\eeqa{1.17a}
and these have the solutions
\beqa \BE_0& = & (\BL_1-\BL_2)^{-1}(\BV_2-\BV_1)\Gt,\quad \BV_*=\BV_1+(\BL_1-\BL_*)(\BL_1-\BL_2)^{-1}(\BV_2-\BV_1),\nonum
c_*& = & fc_1+(1-f)c_2+[f\BV_1+(1-f)\BV_2-\BV_*]^T(\BL_1-\BL_2)^{-1}(\BV_2-\BV_1). \nonum &~&
\eeqa{1.18}
So $c_*$ and $\BV_*$ are determined entirely in terms of $\BL_*$, $f$, and the elements of $\BK_1$ and $\BK_2$. Conversely, if we know
$\BK_*$, then from \eq{1.17} we know $\BL_*$, $\BV_*$ and $c_*$, and the last equation in \eq{1.18} allows us to determine $f$. Thus solving problem 2
is equivalent to solving this problem.

One is often concerned with the quadratic form associated with $\BK_*$ that sometimes may correspond to the energy stored or dissipated in the material. For constant
fields $\BE_0$ and $\Gt$ (with $\BE_0$ not restricted to be given by \eq{1.18}) standard variational principles \cite{Hill:1952:EBC} show that
\beq \bpm \BE_0 \\ \Gt \epm\cdot\BK_*\bpm \BE_0 \\ \Gt \epm=\inf_{\Bu}\bpm \BE_0+\Grad\Bu \\ \Gt \epm\cdot\BK(\Bx)\bpm \BE_0+\Grad\Bu \\ \Gt \epm.
\eeq{1.19}
If we are interested in the lowest value of this over all $\BK_*\in G(\BK_1,\BK_2)$, normalized with say $\Gt=1$, and use an idea of Kohn \cite{Kohn:1991:RDW}, we get
\beqa &~& \inf_{\BK_*\in G(\BK_1,\BK_2)} \bpm \BE_0 \\ 1 \epm\cdot\BK_*\bpm \BE_0 \\ 1 \epm \nonum
& ~ & \quad =\inf_{\Gc}\inf_{\Bu}\lang\bpm \BE_0+\Grad\Bu \\ 1 \epm\cdot[\Gc(\Bx)\BK_1+(1-\Gc(\Bx))\BK_2]\bpm \BE_0+\Grad\Bu \\ 1 \epm\rang \nonum
& ~ & \quad =\inf_{\Bu}\lang\inf_{\Gc}\bpm \BE_0+\Grad\Bu \\ 1 \epm\cdot[\Gc(\Bx)\BK_1+(1-\Gc(\Bx))\BK_2]\bpm \BE_0+\Grad\Bu \\ 1 \epm \rang \nonum
& ~ & \quad =\inf_{\Bu}\lang \min_{j=1,2}\bpm \BE_0+\Grad\Bu \\ 1 \epm\cdot\BK_j\bpm \BE_0+\Grad\Bu \\ 1 \epm\rang \nonum
& ~ & \quad =\inf_{\Bu}\lang W(\BE_0+\Grad\Bu)\rang,
\eeqa{1.20}
where $W(\BF)$ is given by \eq{1.1} and \eq{1.1a}. So we arrive back at the quasiconvexification of $W(\BF)$ as in problem 1, with $m=n$. If $\Gc$ is restricted
to multiple rank laminate geometries we arrive back at the rank-one convexification of $W(\BF)$ (see \cite{Allaire:1999:MDW} and section 31.6 of \cite{Milton:2002:TOC}).
So problem 1 is solved according to whether or not
\beq \inf_{\BK_*\in G(\BK_1,\BK_2)} \bpm \BE_0 \\ 1 \epm\cdot\BK_*\bpm \BE_0 \\ 1 \epm =\inf_{\BK_*\in G^L(\BK_1,\BK_2)} \bpm \BE_0 \\ 1 \epm\cdot\BK_*\bpm \BE_0 \\ 1 \epm.
\eeq{1.21}
\newline
To have equality it is sufficient, but not necessary, to have $G(\BK_1,\BK_2)=G^L(\BK_1,\BK_2)$.

On the other hand, we know the sets $G(\BL_1,\BL_2)$ and $G_f(\BL_1,\BL_2)$ have sufficient convexity (as guaranteed by their stability under lamination) to
be completely characterized by their ``W-transforms''. These generalize the idea of the Legendre transform for characterizing convex sets. First note
that a linear operator $\BA$ on $\CS_{n}$ has elements $A_{ijk\ell}$ such that if the matrix $\BC\in\CS_{n}$ has elements $C_{k\ell}$ then
$\BA\BC$ has elements
\beq \{\BA\BC\}_{ij}=\sum_{k=1}^2\sum_{\ell=1}^n A_{ijk\ell}C_{k\ell}. \eeq{1.21a}
Introducing the inner product
\beq \BA:\BB=\sum_{i,k=1}^2\sum_{j,\ell=1}^nA_{ijk\ell}B_{ijk\ell}, \eeq{1.21b}
between two linear operators $\BA$ and $\BB$ on  $\CS_{n}$, the W-transform of $G(\BL_1,\BL_2)$ is
\beq W(\BN,\BN_\perp)=\inf_{\BL_*\in G(\BL_1,\BL_2)}\{\BN:\BL_*+\BN_\perp:\BL_*^{-1}\},
\eeq{1.22}
where $\BN$ and $\BN_\perp$ range over all real, positive semidefinite, and symmetric operators such that $\BN\BN_\perp=\BN_\perp\BN=0$. When $\BN_\perp=0$ and
$\BN$ is not restricted to be positive semidefinite, this is just
the standard Legendre transform. That  $G(\BL_1,\BL_2)$ may be
characterized in this way is suggested by results of Cherkaev and Gibiansky \cite{Cherkaev:1992:ECB,Cherkaev:1993:CEB} for particular examples
and proved, in general, in \cite{Francfort:1994:SCE}
(see also \cite{Milton:2016:PEE} and section 30.3 of \cite{Milton:2002:TOC}, and references therein).  Writing
\beq \BN=\sum_{k=1}^h\BE_k\otimes\BE_k,\quad  \BN_\perp=\sum_{k=h+1}^n\BJ_k\otimes\BJ_k,
\eeq{1.23}
where some of the $\BE_k$ or $\BJ_k$ could be zero and, without loss of generality, assuming
\beq\BE_k\cdot\BE_\ell=0,\quad\BJ_k\cdot\BJ_\ell=0,\quad \BJ_k\cdot\BE_\ell=0, \text{ for all }k\ne\ell, \eeq{1.23a}
we obtain
\beq \BN:\BL_*+\BN_\perp:\BL_*^{-1}=\sum_{k=1}^h\BE_k\cdot\BL_*\BE_k+\sum_{k=h+1}^n\BJ_k\cdot\BL_*^{-1}\BJ_k.
\eeq{1.24}
Each of the terms in the first sum can be expressed in variational form, similar to \eq{1.19},
\beq \BE_k\cdot\BL_*\BE_k=\inf_{\Bu_k}\lang[\BE_k+\Grad\Bu_k]\cdot\BL(\Bx)[\BE_k+\Grad\Bu_k]\rang,
\eeq{1.25}
while the remaining terms in the second sum can be expressed in the dual
variational form,
\beqa \BJ_k\cdot\BL_*^{-1}\BJ_k& = & \inf_{\Bv_k}\lang[\BJ_k+\BR_\perp\Grad\Bv_k]\cdot[\BL(\Bx)]^{-1}[\BJ_k+\BR_\perp\Grad\Bv_k]\rang\nonum
& = & \inf_{\Bv_k}\lang[\BR_\perp^T\BJ_k+\Grad\Bv_k]\cdot[\BR_\perp\BL(\Bx)\BR_\perp^T]^{-1}[\BR_\perp^T\BJ_k+\Grad\Bv_k]\rang, \nonum &~&
\eeqa{1.26}
where the infimum is over all periodic functions $\Bv_k$, and
\beq
\BR_\perp= \bpm 0 & -1\\ 1 & 0\epm \eeq{1.27}
is the matrix for a $90^\circ$ rotation. 
Let us introduce a constant superfield $\underline{\BE}_0$ and supertensors $\underline{\BL}_1$ and $\underline{\BL}_2$ given by
\beq \underline{\BE}_0=\bpm \BE_1 \\ \vdots \\ \BE_h \\ \BR_\perp^T\BJ_{h+1} \\ \vdots \\ \BR_\perp^T\BJ_n \epm,\quad
\underline{\BL}_j=\bpm \BL_j & \ldots  & 0 & 0 & \ldots & 0 \\
\vdots & \ddots  & \vdots & \vdots & \ddots & \vdots \\
0  & \ldots  & \BL_j & 0 & \ldots & 0 \\
0  & \ldots  & 0 & [\BR_\perp\BL_j\BR_\perp^T]^{-1} & \ldots & 0 \\
\vdots & \ddots  & \vdots & \vdots & \ddots & \vdots \\
0  & \ldots  & 0 & 0 & \ldots & [\BR_\perp\BL_j\BR_\perp^T]^{-1}
\epm.
\eeq{1.28}
Then \eq{1.25} and \eq{1.26} imply
\beq W(\BN,\BN_\perp)=\inf_{\underline{\Bu}}\lang\min_{j=1,2}\{(\underline{\BE}_0+\Grad\underline{\Bu})\cdot\underline{\BL}_j(\underline{\BE}_0+\Grad\underline{\Bu})\}\rang,
\eeq{1.29}
in which the infimum is over all periodic potentials
\beq \underline{\Bu}=\bpm \Bu_1 \\ \vdots \\ \Bu_h \\ \Bv_{h+1} \\ \vdots \\ \Bv_n \epm. \eeq{1.29a}
Thus finding $G(\BL_1,\BL_2)$ is reduced to a set of two-well quasiconvexification problems, each indexed by the value of $h=0,1,\ldots,n$ and with $m=n^2$. 
The problem of finding $G_f(\BL_1,\BL_2)$ can be handled in a similar way \cite{Francfort:1994:SCE}. Instead of \eq{1.22} one considers
\beq W(\BN,\BN_\perp,c)=\inf_f\inf_{\BL_*\in G_f(\BL_1,\BL_2)}\{\BN:\BL_*+\BN_\perp:\BL_*^{-1}+cf\},
\eeq{1.29b}
where the constant $c$ acts as a Lagrange multiplier for the volume fraction $f=\lang\Gc\rang$. One easily sees that this again
reduces to a two-well quasiconvexification problem. 
\newline
\medskip

\section{Some open problems related to the effective conductivity as a function of the component conductivities}
\setcounter{equation}{0}
The lamination closure and the $G_f$-closure coincide when the block entries of $\BL_1$ and $\BL_2$ are all proportional to the $2\times 2$ identity matrix $\BI$,
\beq \BGs^{(k\ell)}_j=\Gs^{(k\ell)}_j\BI,\quad j=1,2. \eeq{1.30}
To see this, we start by following Straley \cite{Straley:1981:TPI} and Milgrom and Shtrikman \cite{Milgrom:1989:LRT} (see also Chapter 6 in \cite{Milton:2002:TOC} and references therein)
and introduce a non-singular matrix $\BW$ having block entries proportional to $\BI$,
\beq \BW=\bpm w^{(11)}\BI & w^{(12)}\BI & \ldots & w^{(1n)}\BI \\
w^{(21)}\BI & w^{(22)}\BI & \ldots & w^{(2k)}\BI & \\
\vdots & \vdots & \ddots & \vdots & \\
w^{(n1)}\BI & w^{(n2)}\BI & \ldots & w^{(nn)}\BI
\epm.
\eeq{1.31}
Now we rewrite \eq{1.10} in the form
\beq \underbrace{\BW^T\BJ(\Bx)}_{\BJ'(\Bx)}=[\Gc(\Bx)\underbrace{\BW^T\BL_1\BW}_{\BL_1'}+(1-\Gc(\Bx)\underbrace{\BW^T\BL_2\BW}_{\BL_2'}]\underbrace{\BW^{-1}\BE(\Bx)}_{\BE'(\Bx)}.
\eeq{1.32}
By choosing $\BW=\BL_2^{-1/2}\BQ$ with $\BQ^T\BQ=\BI$ we get $\BL'_2=\BI$, and then $\BQ$ can be chosen so  $\BL_1'=\BQ^T\BL_2^{-1/2}\BL_1\BL_2^{-1/2}\BQ$ is diagonal,
of the form
\beq \BL_1'=\bpm \Gs^{(1)}\BI & 0 & \ldots & 0 \\
0 & \Gs^{(2)}\BI & \ldots & 0 & \\
\vdots & \vdots & \ddots & \vdots & \\
0 & 0 & \ldots & \Gs^{(n)}\BI.
\epm.
\eeq{1.33}
Thus we have reduced the problem down to a set of uncoupled conductivity problems and the associated effective tensor $\BL'_*=\BW^T\BL_*\BW$ is given by
\beq \BL_*'=\bpm \BGs_*(\Gs^{(1)}) & 0 & \ldots & 0 \\
0 & \BGs_*(\Gs^{(2)}) & \ldots & 0 & \\
\vdots & \vdots & \ddots & \vdots & \\
0 & 0 & \ldots & \BGs_*(\Gs^{(n)})
\epm,
\eeq{1.34}
where $\BGs_*(\Gs)$ is the effective conductivity tensor associated with the equations
\beq \Bj(\Bx)=[\Gc(\Bx)\Gs+(1-\Gc(\Bx))]\Be(\Bx),\quad\Div\Bj=0,\quad \Be=\Be_0+\Grad V, \eeq{1.35}
in which $V(\Bx)$ is a periodic potential, and
\beq \lang\Bj\rang=\BGs_*(\Gs)\lang\Be\rang \eeq{1.36}
defines the function $\BGs_*(\Gs)$. Allowing for complex values of $\Gs$, the properties of this function have been studied in \cite{Bergman:1978:DCC,Milton:1981:BCP,Golden:1983:BEP}.
We remark that complex values of $\Gs$ and hence $\BGs_*$ or, equivalently, complex values of the dielectric constants of the phases and
hence the effective dielectric constant $\BGe_*$ have a physical significance for elecromagnetic waves propagating through the structure when the wavelengths and attenuation
lengths of the waves in each phase are much larger than the microstructure. This is called the quasistatic regime.
In particular $\Imag\BGe_*$ is related to the energy absorption in the composite, and hence is
positive semidefinite when the dielectric constants of the phases are non-negative. Reflecting this,
the function $\BGs_*(\Gs)$ satisfies the Nevanlinna-Herglotz type property, 
\beq \Imag{\BGs_*(\Gs)}\geq 0~ \text{ when}~ \Imag\Gs > 0.
\eeq{1.37}
Additionally, the function is analytic in $\Gs$ except along the negative real $\Gs$-axis,
satisfies the constraints that
\beq \BGs_*(1)=1,\quad \frac{d\BGs_*(\Gs)}{d\Gs}{\Big{|}}_{\Gs=1}=f\BI,\quad \BGs_*(\Gs)\geq 0~ \text{ when}~\Gs~ \text{is real and positive}, \eeq{1.37a}
and, in two-dimensions, the Keller-Dykhne-Mendelson relationship \cite{Keller:1964:TCC, Dykhne:1970:CTD, Mendelson:1975:TEC}, 
\beq \BGs_*(1/\Gs)=\BR_\perp[\BGs_*(\Gs)]^{-1}\BR_\perp^T, \eeq{1.38}
where $\BR_\perp$, with transpose $\BR_\perp^T$ is the matrix for a $90^\circ$ rotation given by \eq{1.27}.
Conversely, any function satisfying these properties can be approximated arbitrarily well by a rational function that corresponds to the effective
conductivity function $\BGs_*^L(\Gs)$ of a hierarchical laminate geometry \cite{Milton:1986:APLG} (see also Section 18.5 in \cite{Milton:2002:TOC}).
Roughly speaking, given this rational function one can retrieve information
about the last two layerings in the corresponding laminate by either setting $\Gs=0$ or $\Gs=\infty$. One strips this last layering away, and accordingly
modifies the associated conductivity function. Then one makes the opposite choice $\Gs=\infty$ or $\Gs=0$, respectively, and proceeds by induction,
until one is left with purely phase 1 or purely phase 2. This establishes that the lamination closure and the $G_f$-closure coincide when the block entries
of $\BL_1$ and $\BL_2$ are all proportional to the $2\times 2$ identity matrix $\BI$. Explicit expressions for the $G_f$-closure were given in
the case $n=1$ by Lurie and Cherkaev \cite{Lurie:1982:AEC} and Murat and Tartar\cite{Murat:1985:CVH} (extended to the three dimensions in
\cite{Lurie:1986:EEC, Murat:1985:CVH}), in the case $n=2$ by Cherkaev and Gibiansky \cite{Cherkaev:1992:ECB},
and for general $n$, using the analytic properties of  $\BGs_*(\Gs)$, by Clark and Milton \cite{Clark:1995:OBC}. It is an open question as to
whether the $G_f$-closure for general $n$ can be obtained via the translation method. One
can speculate that there should be some sort of inductive procedure using the translation method, but it is difficult to see how to formulate this.

In three-dimensions one would like to address the analogous question, and focusing on isotropic composites this becomes:
\newline
\medskip

\noindent
{\it Problem 3: For three-dimensional isotropic composites each having an effective conductivity $\Gs_*\BI$ and being
  built from two isotropic materials having conductivities $\Gs\BI$ and $\BI$,
  can one characterize all possible conductivity functions $\Gs_*(\Gs)$?}
\newline
\medskip

\noindent
The conductivity function $\BGs_*(\Gs)=\Gs_*(\Gs)\BI$ still satisfies \eq{1.37} and \eq{1.37a}, but in place of \eq{1.38} it has been established \cite{Brown:1955:SMP,Bergman:1978:DCC} that
\beq \frac{d^2\Gs_*(\Gs)}{d\Gs^2}{\Big{|}}_{\Gs=1}=-2f(1-f)/3,
\eeq{1.49a}
and, additionally
\cite{Milton:1981:BCP,Avellaneda:1988:ECP,Nesi:1991:MII,Zhikov:1991:EHM,Zhikov:1992:EEH}, that the inequality
\beq \Gs_*(\Gs)\Gs_*(1/\Gs)+\frac{\Gs_*(\Gs)+\Gs\Gs_*(1/\Gs)}{\Gs+1} \geq 2
\eeq{1.50}
holds for all real positive $\Gs$ (and is satisfied as an equality for multicoated sphere assemblages).
The question is whether there exist additional constraints satisfied by $\Gs_*(\Gs)$, and, if so, to identify them. An associated problem is:
\newline
\medskip

\noindent
{\it Problem 4: For three-dimensional isotropic composites of two isotropic phases, are all possible conductivity functions $\Gs_*(\Gs)$ achievable by multiple rank laminate microstructures
  and, if so, does it suffice to consider laminate microstructures where one laminates only in mutually orthogonal directions?}
\newline
\medskip

\noindent
We remark that it does not suffice (even in two-dimensions) to consider laminate microstructures where one laminates in mutually orthogonal directions if one considers anisotropic composites
of two isotropic phases since if $\Gs$ is complex the real and imaginary parts of $\BGs_*(\Gs)$ do not necessarily commute, while they do commute if  one laminates in mutually orthogonal directions.

These results motivate one to consider periodic composites of two anisotropic phases where the conductivity tensor takes the form
\beq \BGs(\Bx)=\Gc(\Bx)\BGs_1+[1-\Gc(\Bx)]\BGs_2, \eeq{1.51}
where the indicator function $\Gc(\Bx)$ is given by \eq{1.12}
and $\BGs_1$ and $\BGs_2$ are the $2\times 2$ matrix-valued conductivity tensors of the two phases. The associated effective
conductivity tensor is found by looking for current fields $\Bj(\Bx)$ and electric fields $\Be(\Bx)$, with the same periodicity of the composite,
that solve
\beq \Bj(\Bx)=\BGs(\Bx)\Be(\Bx),\quad \Div\Bj=0,\quad \Be=-\Grad V(\Bx). \eeq{1.52}
In these equations $V(\Bx)$ is the electric potential, and the volume average, $\lang\Be\rang$, of the electric field $\Be(\Bx)$ is prescribed. Here and
later the angular brackets $\lang\cdot\rang$ denote an average over the unit cell. The average current field $\lang\Bj\rang$ depends linearly
on $\lang\Be\rang$, and it is this linear relation,
\beq \lang\Bj\rang=\BGs_*\lang\Be\rang \eeq{1.53}
that determines the effective tensor $\BGs_*$. We arrive at problem 5, again closely related to problems 1 and 2:
\newline
\medskip

\noindent
{\it Problem 5: For two-dimensional anisotropic composites of two anisotropic phases, are all possible conductivity functions $\BGs_*(\BGs_1,\BGs_2)$ achievable by multiple rank laminate microstructures?}
\newline
\medskip

\noindent
Some progress in characterizing the possible conductivity functions $\BGs_*(\BGs_1,\BGs_2)$ has been made by finding suitable representations of the underlying operators so that they
  satisfy the required algebraic properties \cite{Milton:2018:AET}. Once one has these representations one can, in principle, determine not only $\BGs_*(\BGs_1,\BGs_2)$ but also $\BL_*(\BL_1,\BL_2)$
  for all real positive definite $\BL_1$ and $\BL_2$ taking the block matrix form \eq{1.13}. Thus if one could show a direct correspondence between the operator representations for an
  arbitrary $\Gc(\Bx)$ and the operator representations for  multiple rank laminate microstructures, one would have resolved problem 2,
  establishing that the $G_f$-closure equals its lamination closure.
  Such a correspondence between operator representations was used in \cite{Clark:1994:MEC, Clark:1997:CFR} to establish that in two dimensions the
  effective conductivity function $\BGs_*(\BGs_0)$ of any polycrystal with conductivity of the form
  \beq \BGs(\Bx)=\BR(\Bx)\BGs_0\BR^T(\Bx),\quad \BR(\Bx)\BR^T(\Bx)=\BI, \eeq{1.54}
  and $\BGs_*$ given by \eq{1.52} and \eq{1.53}, corresponds to the conductivity function of a laminate microstructure.

  A question of obvious importance is to identify those two-phase microstructures
  that absorb as much electromagnetic energy as possible, no matter what the
  direction of the incident radiation. In the quasistatic limit, where the
  wavelength of the radiation is much larger than the size of the unit cell of
  periodicity, the electromagnetic equations decouple into separate electric
  equations and magnetic equations involving complex fields and complex
  electrical permittivities and complex magnetic permeabilities, respectively.
  Each decoupled equation is equivalent to a conductivity equation,
  with complex conductivities. Four decades ago bounds were derived on the effective
  complex electrical permittivity (or equivalently the complex magnetic
  permeability, or complex conductivity) of isotropic composites of two
  isotropic phases, mixed in fixed proportions \cite{Bergman:1980:ESM, Milton:1980:BCD}. The bounds confine
  the effective electrical permittivity to a lens-shaped region of the
  complex plane bounded by two circular arcs. The problem becomes one of
  identifying microstructures that have have the maximum imaginary part
  of the effective complex electrical permittivity. In two-dimensions these are assemblages
  of doubly coated disks (corresponding to the transverse electrical
  permittivity of doubly coated cylinders) as they attain the bounds \cite{Milton:1981:BCP}.
  In three-dimensions new bounds \cite{Kern:2020:RCE} show that assemblages of doubly coated
  spheres provide one bounding circular arc. The previously known second bounding arc \cite{Bergman:1980:ESM, Milton:1980:BCD} corresponds to
  conductivity functions $\Gs_*(\Gs)$ that have just one pole at a finite negative real value of $\Gs$.
  Originally just five microgeometries were identified that correspond to five points
  on the circular arc \cite{Milton:1981:BCP}. Depending on the material moduli, these can have
  the maximum possible absorption. Now an extra 3 additional multiple rank laminate geometries have been identified with effective
  electrical permittivities lying on the arc, and which can have the maximum
  possible absorption \cite{Kern:2020:RCE}.  This leads to the following question:
\newline
\medskip

\noindent
  {\it Problem 6: Are there other geometries with isotropic effective permittivities that lie on the arc?}
\newline
\medskip

\noindent
  There is also a close connection with finding isotropic geometries that attain
  bounds on the complex effective bulk modulus \cite{Gibiansky:1993:EVM}, and which can provide
  the maximum possible absorption under oscillatory hydrostatic loadings,
  and that attain bounds coupling the real effective moduli of two conductivity
  type problems that may separately correspond to say, magnetic, thermal,
  particle diffusion, or fluid permeability problems \cite{Bergman:1976:VBS, Bergman:1978:DCC}

  Another question is the following one:
\newline
\medskip

\noindent
  {\it Problem 7: Can any of these discovered geometries, having maximum
  absorption, can be replaced by simpler ones?}
\newline
\medskip

\noindent
  In particular, can the
  assemblages of doubly coated disks or coated spheres be replaced by
  periodic ones with only one inclusion per unit cell? In the case of
  assemblages of coated spheres (isotropic composites having the minimum
  and maximum conductivities for given real positive conductivities of the two
  phases, mixed in given proportions) equivalent periodic geometries having
  only one inclusion per unit cell are known \cite{Vigdergauz:1994:TDG, Grabovsky:1995:MMEb, Liu:2007:PIM}. 
  
\section{Bounds on the elastic moduli of an elastic material with voids, and
  the ultimate auxetic material in this class of materials}
\setcounter{equation}{0}

  Characterizing the possible elasticity tensors of anisotropic composites is a
  daunting task. Elasticity tensors have 18 invariants in three dimensional space and 5 invariants in
  two dimensions, and correspondingly the set of all possible elasticity tensors built from
  two isotropic phases in prescribed volume fractions is represented by a set in an 18 or 5, dimensional
  space, or 21 and 9 if we include the bulk and shear moduli of both phases. The difficulty of this
  is indicated by the observation that a distorted hypercube in 18 dimensions has $2^{18}\approx 26,000$
  vertices and 18 numbers are needed to specify the coordinates of each, bringing the total to about $4.7$
  million numbers, just to specify an 18-dimensional distorted cube. The $G$-closure has only been completely
  characterized, and consists of all positive definite elasticity tensors, in the limit as one phase becomes
  arbitrarily compliant while the other phase  becomes arbitrarily stiff \cite{Milton:1995:WET}. A lot of progress has been made in the
  case where one phase is void, while the other is isotropic has with fixed positive elastic moduli, \cite{Milton:2016:PEE}
  (or when a rigid material replaces the void phase \cite{Milton:2016:TCC}). Still, there are still parts
  of the $G$-closure that have not been mapped. We arrive at
\newline
\medskip

\noindent
  {\it Problem 8: Can one complete the characterization of the $G$-closure for a void (or rigid) phase mixed with an isotropic
  elastic phase?}
\newline
\medskip

\noindent
  It may be the case that the necessary insight for progressing further, at
  least in the case that one phase is void, comes from a consideration
  of the possible pairs of the effective bulk modulus, $\Gk_*$, and effective shear modulus $\Gm_*$, of isotropic composites of an elastic material,
  having bulk and shear moduli $\Gk$ and $\Gm$,  and void. One has the elementary bounds \cite{Hill:1952:EBC}:
  \beq 0\leq \Gk_* \leq \Gk,\quad 0\leq \Gm_* \leq \Gm. \eeq{eb1}

  Naturally the void has minimum effective bulk and shear moduli, both being zero, and the pure elastic phase has
  maximum effective bulk and shear moduli. Also one can construct composites with $(\Gk_*,\Gm_*)$ arbitrarily close
  to  $(\Gk_*,0)$ for all positive $\Gk_*<\Gk$, and  arbitrarily close
  to  $(\Gk,\Gm_*)$  for all positive $\Gk_*<\Gm$ \cite{Milton:2016:PEE, Ostanin:2017:PCC, Milton:2018:NOP} On the other hand, the question remains as
  to what microstructures have high effective shear modulus and low effective bulk modulus. We are led to 
 \newline
\medskip

\noindent
  {\it Problem 9: The bounds \eq{eb1} imply $\Gm_*-c\Gk_*\leq \Gm$ for all $c>0$. Can this inequality be improved, in 2 and/or 3 dimensions,
  for a range of $c>0$? Alternatively, can one construct composites of an elastic phase with void with $(\Gk_*,\Gm_*)$ arbitrarily close to  $(0,\Gm)$?}
\newline
\medskip

\noindent
  A related question is
\newline
\medskip

\noindent
  {\it Problem 10: Identify, for given $c>0$, in 2 and/or 3 dimensions, isotropic microstructures  of an elastic phase with void that have the largest possible value of  $\Gm_*-c\Gk_*$ (or
  a sequence of  isotropic microstructures with moduli such that  $\Gm_*-c\Gk_*$ converges to its largest possible value).}
\newline
\medskip

\noindent
  When $c$ is extremely large, this amounts to identifying  isotropic microstructures that have the largest possible value
  of $\Gm_*$ subject to the constraint that $\Gk_*$ is arbitrarily close to zero. This is what one may call the ultimate auxetic material
  within the class of materials built from an isotropic elastic phase with voids. Auxetic composites have a negative Poisson's ratio, so that they fatten
  when they are pulled, corresponding to a ratio $\Gk_*/\Gm_*<2/3$. When one seeks materials built from an isotropic elastic phase with void, that have
  Poisson's ratios close to the limiting value of $-1$ and thus with $\Gk_*/\Gm_*$ close to zero, it is generally the case that both $\Gk_*$ and $\Gm_*$
  are very small, not just  $\Gk_*$. This is a feature of auxetic composites built from rotating elements \cite{Milton:1992:CMP, Prall:1996:PCH, Grima:2000:ABR} and is less than ideal as one
  wants to retain shear stiffness.
  
   In two-dimensions one can construct a candidate for the title of the ultimate auxetic material as follows. One first takes the elastic phase and slices it into slabs of
  uniform thickness with the interfaces perpendicular to the $x_1$-axis. The slabs are separated by microstructured layers,
  very thin compared to the slab thickness. The microstructured layers are such that their only easy mode of deformation is compression
  of the layer in the direction $x_1$. The thin microstructured layers may, for example, contain the third rank laminate material with a herringbone structure
  depicted in figure 13 of \cite{Milton:1992:CMP} or in the second subfigure of figure 8 in \cite{Milton:2016:PEE}.
  The macroscopic constitutive relation of the sliced material separated by these
  microstructured layers, is 
  \beq \bpm \Gs_{11} \\ \Gs_{22} \\  \Gs_{12} \epm
  = \bpm c_{1111} & c_{1122} & 0\\
  c_{1122} & c_{2222} & 0 \\
  0 & 0 & 2c_{1212} 
  \epm \bpm  \Ge_{11} \\ \Ge_{22} \\ \Ge_{12} \epm,
  \eeq{aux1}
  where the $\Gs_{ij}$ are the Cartesian components of the average stress, while the $\Ge_{ij}$ are the Cartesian components of the average strain. 
  The effective elastic moduli are 
  \beq c_{1111}=\Gve,\quad  c_{1122}=c\Gve,\quad c_{2222}\approx\frac{4\Gk\Gm}{\Gk+\Gm}=E,\quad c_{1212}\approx \Gm,
  \eeq{aux2}
  where $\Gve$ is a small parameter, reflecting the easyness of the easy mode of compression in the $x_1$-direction,
  and the appearance of $E=4\Gk\Gm/(\Gk+\Gm)$ reflects the fact that the effective Young's modulus for compression in
  the $x_2$-direction is approximately the same as the pure elastic phase, namely $E$. 
  We now treat this material as a crystal and construct from it the polycrystal with the largest possible effective shear modulus $\Gm_*$ and
  smallest possible effective bulk modulus $\Gk_*$. According to the bounds and laminate constructions in \cite{Avellaneda:1996:CCP} these are
  \beqa \Gk_*& = & \frac{c_{1111}c_{2222}-c_{1212}^2}{c_{1111}+c_{2222}-2c_{1212}}, \nonum
  \quad \Gm_*& = &\frac{c_{1111}c_{2222}-c_{1212}^2}{2c_{1212}-2c_{2222}+2\sqrt{c_{2222}[c_{1111}+c_{2222}-2c_{1212}+(c_{1111}c_{2222}-c_{1212}^2)/c_{1212}]}}. \nonum
  &~&
  \eeqa{aux3}
  Substituting \eq{aux2} in these, and taking the limit $\Gve\to 0$ gives
  \beq \Gk_*=0, \quad \frac{1}{\Gm_*}=\frac{5}{4\Gm}+\frac{1}{4\Gk}. \eeq{aux4}
  The formula for $\Gm_*$ has the required invariance property that if $1/\Gm$ and $-1/\Gk$ are shifted by the same constant, then
  $1/\Gm_*$ is shifted by this constant too \cite{Lurie:1984:CSP, Cherkaev:1992:IPS}. Due to this invariance we may assume, without loss of generality, that the initial elastic
  phase is incompressible ($1/\Gk=0$) so that \eq{aux4} implies $\Gm_*=4\Gm/5$. The question is then:
\newline
\medskip

\noindent
  {\it Problem 11: Is $4\Gm/5$ the largest  possible value of $\Gm_*$ for a two-dimensional elastic material with voids, given that $\Gk_*=1/\Gk=0$?}
\newline
\medskip

\noindent
  From a practical standpoint the answer to this question is moot, as not only are such multiple rank laminates impossible to build and
  subject to buckling, but also the linear elastic moduli are largely irrelevant under finite but small deformations as the microstructured layers
  will undergo large deformations relative to their thickness. Ideally one wants to address
\newline
\medskip

\noindent
  {\it Problem 12: Can one obtain bounds that correlate the possible compressive and shear deformations of composites when these deformations are not infinitesimal?}  
  \newline
\medskip

\noindent
Returning back to the theoretical problem of finding the ultimate auxetic material, one could use
  in principle a similar construction in three-dimensions. 
  However the barrier is that the polycrystals having the largest $\Gm_*$ with $\Gk_*=0$ have not yet been identified. Thus one arrives at
\newline
\medskip

\noindent
  {\it Problem 13: What are the possible $(\Gk_*,\Gm_*)$-pairs for three-dimensional isotropic elastic polycrystals (composites built from a single crystal in various
    orientations)? The bounds of Hill \cite{Hill:1952:EBC} are optimal for $\Gk_*$ \cite{Avellaneda:1989:OBE}, but improved bounds for $\Gm_*$ or $(\Gk_*,\Gm_*)$ pairs
    are lacking.  Hashin and Shtrikman obtained improved bounds on  $\Gm_*$ \cite{Hashin:1962:VATE}, but only under additional assumptions about  crystal orientations, that are
  not generally valid.}
  \newline
\medskip

\noindent
For conductivity the analogous problem has been solved \cite{Schulgasser:1977:BCS, Avellaneda:1988:ECP, Nesi:1991:PCM}, but the $G$-closure containing all possible effective
conductivity tensors of anisotropic polycrystals has not yet been fully mapped. 

   More generally, moving back to isotropic composites of two isotropic elastic phases, one possibly rigid or void, the tightest known bounds
  on the possible $(\Gk_*,\Gm_*)$-pairs
  are those of Cherkaev and Gibiansky \cite{Cherkaev:1993:CEB}, in two-dimensions, and those of Berryman and Milton \cite{Berryman:1988:MRC}, in three-dimensions. Its seems highly likely that these
  bounds are not optimal. Gal Shmuel and myself are progressing on a nontrivial route for improving the three dimensional bounds using the 
  ``translation method'' approach (see Chapters 24 and 25 of \cite{Milton:2002:TOC} and references therein)
  used by Cherkaev and Gibiansky, but even so these improved bounds are unlikely to be optimal. Thus we come to
  \newline
\medskip

\noindent
  {\it Problem 14: Can one obtain improved bounds on the elastic moduli pairs of isotropic composites of two isotropic elastic phases, and ultimately
    find the optimal ones?}
 \newline
 \medskip
 
  Numerical explorations of the possible $(\Gk_*,\Gm_*)$-pairs have been made, for example in \cite{Andreassen:2014:DME, Ostanin:2017:PCC}. From a practical
  viewpoint such numerical explorations are probably more useful than the theoretical developments. On the other hand, it is difficult to
  numerically explore multiscale structures that may be necessary to obtain desired extreme responses, such as in resolving Problem 11.

\section{Some future directions for wave and other equations}
\setcounter{equation}{0}

  An impressive body of research addresses the problem of bounding the response of bodies to electromagnetic
  or other waves, and addressing limitations to how one can manipulate these waves. A few examples include the
  results in \cite{Sohl:2007:PLB, Gustafsson:2007:PLA, Rechtsman:2009:MOU, Miller:2015:FLO, Shim:2019:MFS} and references therein.
  There are many problems to be addressed and new approaches are needed to improve existing bounds, or to reveal novel ones.
  A framework suited to most linear equations in physics \cite{Milton:2020:UPLI, Milton:2020:UPLII, Milton:2020:UPLIII, Milton:2020:UPLIV}, including wave and diffusion equations,
  is to express them in the
  form
  \beq \BJ(\Bx)=\BL(\Bx)\BE(\Bx)-\Bs(\Bx),\quad \BJ\in\CJ,\quad\BE\in\CE, \eeq{form}
  where the first equation is the constitutive law, with the tensor $\BL(\Bx)$ representing the local material properties,
  $\Bs(\Bx)$ is the source term, while $\CE$ and $\CJ$ are orthogonal spaces embodying the differential constraints on the
  fields. Here $\Bx$ represents a point in space, or space time with $x_0$ representing time. Scattering problems
  can also be expressed in this form \cite{Milton:2017:BCP} by incorporating the fields ``at infinity'' appropriately.
  The analog for quadratic forms of quasiconvexity is then $Q^*$-convexity: a quadratic
  form $f(\BP)$ is  $Q^*$-convex if $f(\BE)\geq 0$ for all $\BE\in\CE$. $Q^*$-convex functions allow one to place bounds
  on the spectrum of the operator relevant to the problem \cite{Milton:2019:NRF, Milton:2020:UPLV}. The subject of $Q^*$-convexity remains to be explored,
  and simple examples of  $Q^*$-convex functions need to be found for the various equations, beyond quasiconvex functions and those discovered for
  the Schr\"odinger equation (Sections 13.6 and 13.7 of \cite{Milton:2016:ETC}) .  For wave and diffusion equations it seems likely that they will provide a powerful tool for
  addressing other bounding problems, and this provides an avenue for future work. In connection
  with this, variational principles have been
  developed for acoustic, elastic, and electromagnetic equations at constant frequency in lossy materials
  \cite{Milton:2009:MVP,Milton:2010:MVP}. These are the direct
  analogs of those of Gibiansky and Cherkaev \cite{Cherkaev:1994:VPC} that have proved very powerful, in conjunction with the use of quasiconvex functions, for obtaining bounds on the
  quasistatic response of
  composites: examples include bounds on effective complex electrical permittivities (Section 22.6 of \cite{Milton:2002:TOC} and \cite{Kern:2020:RCE})
  and bounds on complex bulk moduli \cite{Gibiansky:1993:EVM}. So one expects there should be useful
  bounds resulting from these variational principles for wave equations in lossy media.

  Recently it has been discovered that associated with exact relations for composites, as reviewed in Chapter 17 of \cite{Milton:2002:TOC} and the book \cite{Grabovsky:2016:CMM},
  are exact relations satisfied by the infinite body Green's function in certain inhomogeneous media, and boundary field equalities \cite{Milton:2019:ERG}.
  Boundary field equalities are exact identities satisfied by the fields at the boundary of the body, given that the fields
  in the interior of the body satisfy some constraints that do not uniquely determine the interior fields in terms
  of their boundary values. A classical example is that a field with zero divergence has zero net flux through the boundary.
  The theory of these exact relations for the Green's function and boundary field equalities extends to wave and diffusion equations \cite{Milton:2019:ERG},
  or more generally to equations expressible in the form \eq{form}, but examples, and in particular useful examples, need to be generated. 

  Another topic to be explored is that of neutral Inclusions for wave equations. For static and quasistatic problems there are many studies of neutral inclusions (see, for example,
  Section 7.11 of \cite{Milton:2002:TOC}, the review \cite{Ji:2020:NIW}, and references therein).  These are inclusions that one
  can insert in a homogeneous medium without disturbing the surrounding fields, provided these fields fall into an appropriate class. Thus, for example, one may
  obtain neutrality for a single applied uniform fields, for any uniform field, or for any applied field satisfying the underlying equations. For conductivity,
  or equivalently for the dielectric problem, coated ellipsoids can be neutral and invisible to any uniform field \cite{Kerker:1975:IB}. In two-dimensions there are other shaped inclusions
  that can be neutral to a uniform field in a specified direction \cite{Milton:2001:NIC}. Coated dielectric cylinders, where the core, coating, and surrounding medium have dielectric
  constants of 1, $-1+i\Gd$, and 1 become neutral and hence invisible to large classes of fields in the limit $\Gd\to 0$ \cite{Nicorovici:1994:ODP}, and can cloak sources and objects
  \cite{Milton:2006:CEA, Nguyen:2017:CAO}. Transformations allow one to obtain other inclusions
  that are neutral and thus invisible to any exterior field, and also cloak objects \cite{Greenleaf:2003:ACC}. The transformation approach also yields neutral inclusions that are invisible to constant frequency
  electromagnetic waves \cite{Dolin:1961:PCT}. Even appropriately coated spheres can be invisible in the far field when the incident is planar \cite{Alu:2005:ATP}.
  Quite simple inclusions have been found that are neutral and hence invisible to a single incident planar electromagnetic wave \cite{Xi:2009:ODP, Landy:2013:FPU}.
  One, possibly difficult, 
  research direction, is to explore whether there are other simple geometries, not obtained from a transformation approach, that are invisible
  to one or more incident plane waves.

  Most analysis of wave equations in lossy media has been done at constant frequency, which makes sense as this avoids convolutions in time. However recent work
  on bounds in the time domain \cite{Mattei:2016:BRV, Mattei:2020:EPA}
  show that it is possible for the temporal response of a two-phase mixture to be untangled at specific times when the applied
  field has an appropriately tailored dependence on time. This shows it may be productive to depart from focusing on bounds at constant frequency, and to consider
  bounding responses as a function of time. Beyond the analytic approach used in these papers, the variational approach of Carini and Mattei \cite{Carini:2015:VFL},
  may be helpful if one can modify it to obtain bounds at each instant in time, rather than to bounding the response over at interval of time.

\section*{Acknowledgment}
The author would like to thank the Isaac Newton Institute for Mathematical Sciences for support and hospitality during the programme ``The Mathematical Design of New Materials``
when work on this paper was initiated. This work was supported by: EPSRC grant number EP/R014604/1. Additionally the author 
is grateful to the National Science Foundation for support through the Research Grant DMS-1814854, and thanks Yury Grabovsky for helpful comments. 

\ifx \bblindex \undefined \def \bblindex #1{} \fi\ifx \bbljournal \undefined
  \def \bbljournal #1{{\em #1}\index{#1@{\em #1}}} \fi\ifx \bblnumber
  \undefined \def \bblnumber #1{{\bf #1}} \fi\ifx \bblvolume \undefined \def
  \bblvolume #1{{\bf #1}} \fi\ifx \noopsort \undefined \def \noopsort #1{}
  \fi\ifx \bblindex \undefined \def \bblindex #1{} \fi\ifx \bbljournal
  \undefined \def \bbljournal #1{{\em #1}\index{#1@{\em #1}}} \fi\ifx
  \bblnumber \undefined \def \bblnumber #1{{\bf #1}} \fi\ifx \bblvolume
  \undefined \def \bblvolume #1{{\bf #1}} \fi\ifx \noopsort \undefined \def
  \noopsort #1{} \fi

\end{document}